\documentclass[11pt,a4paper]{article}
\usepackage[margin=1.0in]{geometry}
\IfFileExists{lmodern.sty}{\usepackage[T1]{fontenc}\usepackage{lmodern}}{\usepackage[OT1]{fontenc}}
\usepackage[utf8]{inputenc}
\usepackage{microtype}
\microtypesetup{expansion=false}
\usepackage{amsmath,amssymb,amsthm,mathtools}
\usepackage{enumitem}
\usepackage{hyperref}
\hypersetup{colorlinks=true,linkcolor=blue!45!black,citecolor=blue!45!black,urlcolor=blue!45!black,hypertexnames=false,pdftitle={Intermediate Growth and Rigidity in Primitive Circulant Spectral Algebras},pdfauthor={Juan D. Velez and Carlos A. Cadavid}}
\usepackage{xcolor}
\allowdisplaybreaks[2]
\newtheorem{theorem}{Theorem}[section]
\newtheorem{lemma}[theorem]{Lemma}
\newtheorem{proposition}[theorem]{Proposition}
\newtheorem{corollary}[theorem]{Corollary}
\theoremstyle{remark}
\newtheorem{remark}[theorem]{Remark}

\newcommand{\Q}{\mathbb Q}
\newcommand{\Z}{\mathbb Z}

\newcommand{\prim}{\mathrm{prim}}
\newcommand{\glob}{\mathrm{glob}}
\newcommand{\ASg}{A_S^{\glob}}
\newcommand{\ASp}{A_S^{\prim}}

\newcommand{\rhoSp}{\rho_S^{\prim}}
\newcommand{\one}{\mathbf 1}
\newcommand{\gr}{\operatorname{gr}}
\newcommand{\Cay}{\operatorname{Cay}}
\newcommand{\lcm}{\operatorname{lcm}}

\newcommand{\rank}{\operatorname{rank}}

\title{\textbf{Intermediate Growth and Rigidity in Primitive Circulant Spectral Algebras}}

\author{
Juan D. V\'elez\\[-2pt]
\small Universidad Nacional de Colombia, Sede Medell\'in, Colombia
\and
Carlos A. Cadavid\\[-2pt]
\small Universidad EAFIT, Medell\'in, Colombia
}

\date{November 2026}
\begin{document}
\maketitle

\begin{abstract}
At prime graph orders, the complete nontrivial and primitive Fourier
spectra of a circulant coincide. We show that this agreement conceals
a strong algebraic separation when the graph order varies.
For any nonempty finite jump set $S$, let $A_S^{\prim}$ be the algebra
of arithmetic functions obtained by evaluating symmetric power sums
on the primitive spectrum, filtered by power-sum weight.
Removing the common Euler-totient term gives a reduced, locally
finite algebra $B_S$, and we prove a filtration-preserving
isomorphism $A_S^{\prim}\cong B_S[T_2]$.
Thus $A_S^{\prim}$ is non-Noetherian and has Krull dimension one.
Its filtered dimensions satisfy
$\exp(c_S\sqrt{D/\log D})\leq \rho_S^{\prim}(D)
\leq\exp(C_S\sqrt D)$ for large $D$, and its associated-graded
Hilbert series is nonrational.
Despite agreement on every prime order, the primitive and full
spectral evaluation algebras intersect only in $\mathbb Q$.
More strongly, no nonconstant primitive spectral function is
eventually quasirational on residue classes.
We give an exact finite-matrix rank formula and determine the
quotient obtained by restricting to integers with large prime factors.
\end{abstract}

\section{Introduction: the arithmetic distinction between two spectra}
Fix a nonempty finite set $S\subset\mathbb Z_{>0}$, and put
$r=|S|$ and $M=\max S$. For $n>2M$ consider the circulant graph
\[
\Gamma_{n,S}=\Cay(\mathbb Z/n\mathbb Z,\{\pm s:s\in S\}),
\qquad
L_S(z)=\frac12\sum_{s\in S}(z^s+z^{-s}),
\qquad \zeta_n=e^{2\pi i/n}.
\]
Its normalized adjacency eigenvalues are $L_S(\zeta_n^k)$.
Two natural ways of collecting them are the complete nontrivial
Fourier spectrum, indexed by $1\le k<n$, and the primitive spectrum,
indexed by the units modulo $n$. Define the power-sum moments
\begin{align}
G_h(n)&=\sum_{k=1}^{n-1}L_S(\zeta_n^k)^h,
\label{eq:G}\\
T_h(n)&=\sum_{k\in(\mathbb Z/n\mathbb Z)^\times}
L_S(\zeta_n^k)^h.\label{eq:T}
\end{align}
The full and primitive moments agree at every prime order $p>2M$,
because every nonzero residue class modulo $p$ is invertible.
It is tempting to expect that their polynomial spectral invariants
remain closely related. Our main result shows that this expectation
fails dramatically once composite orders are taken into account.

Let $\Lambda_{\mathbb Q}=\mathbb Q[p_1,p_2,\ldots]$,
$\deg p_h=h$. Evaluating the power sums at these two sets of
moments gives filtered algebras of rational arithmetic functions,
\[
A_S^{\glob}=\mathbb Q[G_1,G_2,\ldots],\qquad
A_S^{\prim}=\mathbb Q[T_1,T_2,\ldots]
\subseteq\mathbb Q^{\{n>2M\}}.
\]
Write $\rho_S^\bullet(D)=\dim_{\mathbb Q}A_{S,\le D}^\bullet$
for the dimension of the image of weight at most $D$.

The difference originates in Fourier orthogonality. The full moment
$G_h(n)$ is affine in $n$ beyond $n>hM$, with only finitely many
lower-order corrections. By contrast, the primitive moment is expressed
through Ramanujan sums: its main term is a multiple of Euler's
totient $\varphi(n)$, and its correction detects the small prime
divisors of $n$. Increasing the weight reveals new independent
divisibility conditions. The contrast between finite winding
corrections and prime-factor sensitivity motivates the algebraic
questions studied here.

\paragraph{A motivating example.}
Already $S=\{1\}$ displays the distinction. Here
$L_S(z)=(z+z^{-1})/2$, and for all admissible $n>2$,
\[
G_1(n)=-1,\qquad T_1(n)=c_n(1)=\mu(n).
\]
At each odd prime $p$, both functions take the value $-1$;
at $p^2$, however, the primitive moment is zero, whereas the
full moment remains $-1$. This elementary discrepancy does not
yet explain how far apart the \emph{algebras} are: polynomial
expressions in higher moments could, in principle, impose unexpected
relations. Determining the size of the resulting filtered algebras
and their common elements is the problem addressed below.

\paragraph{Why this filtration?}
The weight $\deg p_h=h$ is the usual power-sum grading of the
ring of symmetric functions~\cite{StanleyEC2}. It measures the
number of factors in a spectral moment and is preserved by
polynomial spectral evaluation. Consequently the filtered dimension
counts independent spectral expressions of bounded combinatorial
complexity, even though the full algebras involve infinitely many
different graph orders.

\begin{theorem}[Universal intermediate filtered growth]\label{thm:growth}
For every nonempty finite $S$, there exist $c_S,C_S>0$ such that
for all sufficiently large $D$,
\begin{equation}\label{eq:prim-bounds}
\exp\left(c_S\sqrt{\frac D{\log D}}\right)
\le \rho_S^{\prim}(D)\le\exp(C_S\sqrt D).
\end{equation}
The associated-graded Hilbert series of $A_S^{\prim}$ is nonrational.
\end{theorem}

There is an accompanying structural description. The functions
obtained by subtracting the common $\varphi(n)$ term from the
primitive moments take only finitely many values; the algebra
they generate is locally finite over $\mathbb Q$. Nevertheless,
adjoining one degree-two moment gives a polynomial extension.
This accounts for Krull dimension one and for the coexistence of
non-Noetherian structure with intermediate filtered growth.
The last phrase requires care: the algebra is not finitely generated,
and its growth refers to this specific power-sum filtration.
It is different from the ordinary word growth used in
Gelfand--Kirillov dimension~\cite{KrauseLenagan}.

For the rigidity statement, let $\mathrm{QR}$ be the algebra of
arithmetic functions that are, eventually on every class modulo
some fixed integer, given by rational functions of $n$ with rational
coefficients.
\begin{theorem}[Arithmetic rigidity]\label{thm:rigidity-intro}
For every nonempty finite $S$,
\begin{equation}\label{eq:intersection}
A_S^{\prim}\cap\mathrm{QR}=\mathbb Q,
\qquad A_S^{\glob}\cap A_S^{\prim}=\mathbb Q.
\end{equation}
In addition, $A_S^{\prim}$ has no nonzero finitely supported element,
although $G_h(p)=T_h(p)$ for all primes $p>2M$ and all $h$.
\end{theorem}

The pointwise equality on primes is therefore insufficient to
determine the ambient arithmetic-function algebra. Prime-square
levels supply a decisive test: sufficiently large prime squares
annihilate all nonzero-frequency Ramanujan corrections at any fixed
weight, and the resulting polynomial identity is constrained by two
incompatible reflections. This yields the rigidity theorem without
assuming an asymptotic estimate for the primitive dimension.

Our treatment extends the cyclotomic evaluation viewpoint of
V\'elez and Cadavid~\cite{VelezCadavid2026,VelezCadavidOrbit}.
We use standard Fourier orthogonality, Ramanujan sums and M\"obius
inversion~\cite{Apostol,Washington}, and the partition-function
bound of Hardy and Ramanujan~\cite{HardyRamanujan}. The central
questions concern entire filtered algebras as the modulus varies,
rather than a spectral computation at one fixed graph order.

Section~2 establishes the moment bridge. Section~3 describes the
full-spectrum algebra for comparison; Sections~4--5 determine the
structure and growth of the primitive algebra. Section~6 proves
the rigidity theorem; Sections~7--8 describe its rough-level
quotient and a finite-rank evaluation algorithm.

\section{Moments, Fourier orbits, and arithmetic decomposition}
Write the Laurent expansion as
\begin{equation}\label{eq:Laurent}
L_S(z)^h=\sum_{|m|\le hM}b_{h,m}z^m,\qquad a_h=b_{h,0}.
\tag{2.1}\end{equation}
All $b_{h,m}$ are nonnegative rational numbers, $b_{h,m}=b_{h,-m}$, and
\[
a_1=0,\qquad a_2=r/2>0,\qquad \sum_m b_{h,m}=r^h.
\]
The next lemma gives the basic relation between the two kinds of moments.
\begin{lemma}[Exact moment bridge]\label{lem:bridge}
For $n>2M$ and $h\ge1$,
\begin{align}
G_h(n)&=a_hn-r^h+n\sum_{\substack{m\ne0\\n\mid m}}b_{h,m},\label{eq:G-bridge}\tag{2.2}\\
T_h(n)&=a_h\varphi(n)+R_h(n),\qquad R_h(n)=\sum_{0<|m|\le hM}b_{h,m}c_n(m),\label{eq:T-bridge}\tag{2.3}\\
G_h(n)+r^h&=\sum_{d\mid n}T_h(d),\label{eq:divisor-bridge}\tag{2.4}
\end{align}
where $c_n(m)=\sum_{u\in(\Z/n\Z)^\times}\zeta_n^{um}$ is the Ramanujan sum. In particular,
$G_h(n)=a_hn-r^h$ whenever $n>hM$.
\end{lemma}
\begin{proof}
Interchange the two finite sums in the definition of $G_h$ and use
$\sum_{k=0}^{n-1}\zeta_n^{km}=n\one_{n\mid m}$. We must subtract the term corresponding to the trivial mode, which is $L_S(1)^h=r^h$. This gives the first formula. Restricting the same computation to unit indices gives the second. Finally, every $n$th root of unity has one definite order $d$ dividing $n$. Collecting terms according to that order proves the third identity.
\end{proof}
Recall that
\[
c_n(m)=\sum_{e\mid(n,m)}e\mu(n/e).
\]
If
\[
\eta_e(n)=\one_{e\mid n}\mu(n/e),\qquad
A_{h,e}=e\sum_{\substack{0<|m|\le hM\\e\mid m}}b_{h,m},
\]
then
\begin{equation}\label{eq:R-eta}
R_h(n)=\sum_{e=1}^{hM}A_{h,e}\eta_e(n).
\tag{2.5}\end{equation}
It follows that each $R_h$ has finite image. We stress one point: the divisor relation~\eqref{eq:divisor-bridge} is a relation between moments. Dirichlet convolution does not respect pointwise multiplication, so it cannot be used to identify the two algebras.

\begin{lemma}[Insertion of large primes]\label{lem:insertion}
Fix $h\ge1$ and $d\ge1$. If $p>hM$ is prime and $p\nmid d$, then
\begin{equation}\label{eq:insert1}
R_h(dp)=-R_h(d),\qquad T_h(dp)=a_h\varphi(d)p-T_h(d).
\tag{2.6}\end{equation}
For two distinct primes $p,q$ of this kind,
\begin{equation}\label{eq:insert2}
R_h(dpq)=R_h(d),\qquad T_h(dpq)=a_h\varphi(d)(p-1)(q-1)+R_h(d).
\tag{2.7}\end{equation}
For any sufficiently large prime $p$ coprime to $d$,
\begin{equation}\label{eq:insert-square}
R_h(dp^2)=0,\qquad T_h(dp^2)=a_h\varphi(d)p(p-1).
\tag{2.8}\end{equation}
\end{lemma}
\begin{proof}
For every nonzero frequency appearing in~\eqref{eq:T-bridge}, the inequality $p>hM$ gives $p\nmid m$. Hence $c_p(m)=-1$ and $c_{p^2}(m)=0$. The formulas now follow from the multiplicativity of Ramanujan sums in coprime moduli and the corresponding formulas for the Euler totient.
\end{proof}
\begin{remark}[Convention for small divisors]
In~\eqref{eq:divisor-bridge} the symbol $T_h(d)$ is defined for all $d\ge1$ by the formula with the unit group. We regard $(\Z/1\Z)^\times$ as the one-element group, and thus $T_h(1)=L_S(1)^h=r^h$. Only the divisor identity uses this extension. The evaluation algebras themselves are defined on $n>2M$.
\end{remark}

\section{The global algebra as a comparison}
The full-spectrum moments have a simple polynomial part. It is convenient to separate it from the terms coming from winding walks. Put
\begin{equation}\label{eq:tau-Q}
\tau_S=\frac2r(G_2+r^2)=n,\qquad
Q_h(n)=G_h(n)-a_hn+r^h
=n\sum_{\substack{m\ne0\\n\mid m}}b_{h,m},\quad h\ge3.
\tag{3.1}\end{equation}
Here $\tau_S=n$ means equality as arithmetic functions on the admissible levels. The correction $Q_h$ is nonnegative and supported in the finite interval $2M<n\le hM$.

\begin{proposition}[The complete global image]\label{prop:glob-image}
For every nonempty finite $S$, the algebra $\ASg$ is exactly the algebra of rational-valued sequences on $n>2M$ that are eventually polynomial in $n$. An element of filtered weight at most $D$ agrees, for all sufficiently large $n$, with a polynomial of degree at most $\lfloor D/2\rfloor$.
\end{proposition}
\begin{proof}
For $n>DM$, each moment $G_h$ with $h\le D$ has the form $a_hn-r^h$. Since $G_1$ has zero slope, every monomial of total weight at most $D$ gives a polynomial of degree at most $\lfloor D/2\rfloor$. Conversely, the function $\tau_S=n$ belongs to our algebra, so all polynomial sequences belong to it.

We also have all finitely supported sequences. Indeed, for $m>2M$ the coefficient $b_{m,mM}=2^{-m}$ is positive, and hence $Q_m(m)>0$. As $Q_m$ has finite support, interpolation gives the indicator of the level $m$:
\[
E_{m,S}(n)=\frac{Q_m(n)}{Q_m(m)}
\prod_{\substack{2M<j\le mM\\j\ne m}}\frac{n-j}{m-j}.
\]
This belongs to $\ASg$. Combining finitely many such indicators with a polynomial produces any eventually polynomial sequence.
\end{proof}

\section{The primitive algebra as a polynomial extension}
For $h\ge3$, remove from $T_h$ its common totient part by setting
\begin{equation}\label{eq:F-def}
F_h=T_h-\frac{a_h}{a_2}T_2=R_h-\frac{a_h}{a_2}R_2,
\qquad B_S=\Q[T_1,F_3,F_4,\ldots].
\tag{4.1}\end{equation}
Each generator of $B_S$ takes only finitely many values, and $F_h$ has weight at most $h$. These simple facts have useful algebraic consequences.

\begin{lemma}[Independence of the reduced traces]\label{lem:F-independent}
The functions $1,F_3,F_4,\ldots$ are linearly independent over $\Q$ on the levels $n>2M$.
\end{lemma}
\begin{proof}
The functions $\eta_e(n)=\one_{e\mid n}\mu(n/e)$ are linearly independent on all positive integers. To see this, convolve a finite linear relation with the constant-one arithmetic function; the coefficient at each $e$ is then recovered directly. In~\eqref{eq:R-eta}, the coefficient of $\eta_{hM}$ in $F_h$ is
\[
hM\bigl(b_{h,hM}+b_{h,-hM}\bigr)=hM\,2^{1-h}\ne0.
\]
No $F_j$ with $j<h$ has this term. Hence the $F_h$ are independent on all positive integers.

It remains to pass to the admissible levels. Suppose
$u_0+\sum_{h=3}^H u_hF_h(n)=0$ for every $n>2M$. Take $d\ge1$ and evaluate this relation at $dp$ and $dpq$, where $p,q$ are distinct and arbitrarily large primes. Lemma~\ref{lem:insertion} gives
\[
u_0-\sum_{h=3}^H u_hF_h(d)=0,
\qquad u_0+\sum_{h=3}^H u_hF_h(d)=0.
\]
Thus $u_0=0$, and the independence on all positive integers finishes the argument.
\end{proof}

\begin{theorem}[Polynomial-extension structure]\label{thm:polynomial-extension}
For every nonempty finite $S$, the homomorphism
\begin{equation}\label{eq:polynomial-extension}
B_S[X]\longrightarrow\ASp,\qquad X\longmapsto T_2,
\tag{4.2}\end{equation}
is an isomorphism of filtered $\Q$-algebras, with $X$ of weight two and $F_h$ of weight $h$. Every finitely generated subalgebra of $B_S$ is a finite product of copies of $\Q$. In particular, $\ASp$ is reduced, non-Noetherian, integral over $\Q[T_2]$, and has Krull dimension one. Its Gelfand--Kirillov dimension, interpreted through finitely generated subalgebras, equals one.
\end{theorem}
\begin{proof}
Take finitely many rational-valued functions with finite images. They have only finitely many joint values in $\Q^t$. Polynomial interpolation produces the indicator of each such value, so the algebra they generate is a finite product of copies of $\Q$. This proves the assertion about $B_S$.

Surjectivity of~\eqref{eq:polynomial-extension} follows immediately from $T_h=F_h+(a_h/a_2)T_2$. For injectivity, assume
\[
\sum_{j=0}^J u_jT_2^j=0,\qquad u_j\in B_S.
\]
Fix an admissible $d>2M$ and choose distinct primes $p,q$ sufficiently large. By the insertion formulas,
\[
u_j(dpq)=u_j(d),\qquad
T_2(dpq)=a_2\varphi(d)(p-1)(q-1)+R_2(d).
\]
Keep $q$ fixed and vary $p$ through infinitely many primes. The relation is then a polynomial in a variable taking infinitely many different values. All its coefficients must vanish, and therefore $u_j(d)=0$ for every $j$. Since $d$ was arbitrary, all $u_j$ are zero.

We spell out the filtration assertion. Give the formal variables
$Y_1,Y_2,Y_3,\ldots$ weights $1,2,3,\ldots$, respectively.
The triangular substitutions
\[
Y_1\longmapsto T_1,\qquad Y_2\longmapsto T_2,
\qquad Y_h\longmapsto F_h+\frac{a_h}{a_2}T_2
\quad(h\ge3)
\]
do not increase weight, since the extra term has weight two.
Conversely, the substitutions
$F_h=T_h-(a_h/a_2)T_2$ do not increase weight.
Thus these changes of generators induce exactly the same image
filtration. Together with injectivity of~\eqref{eq:polynomial-extension},
this proves the filtered polynomial-extension assertion.

All elements of $B_S$ are integral over $\Q$, while $T_2$ is
transcendental: it already has infinitely many values on prime
levels. Integral extensions preserve Krull dimension, so the
dimension is one~\cite{Stacks}.

By Lemma~\ref{lem:F-independent}, the polynomial extension contains infinitely many elements linearly independent over $\Q[T_2]$, and so it is not module-finite over this ring. To see that it is non-Noetherian, observe first that $B_S$ is reduced and integral over $\Q$, hence zero-dimensional. If $B_S$ were Noetherian, it would be Artinian and therefore a finite product of fields. Each field factor embeds, by evaluation at a suitable admissible level, into $\Q$; thus every factor would equal $\Q$. This would make $B_S$ finite-dimensional, contrary to Lemma~\ref{lem:F-independent}. Hence $B_S$ is not Noetherian, and neither is $B_S[T_2]$, since its quotient by $(T_2)$ is $B_S$. Reducedness of the polynomial extension follows because it is an algebra of functions. Finally, every finitely generated subalgebra is contained in $C[T_2]$ for some finite-dimensional $C\subset B_S$, so its word growth is at most linear. The subalgebra $\Q[T_2]$ gives the opposite bound.
\end{proof}

Let $b_S(D)=\dim_{\Q}(B_S)_{\le D}$, with $b_S(D)=0$
for $D<0$; here $B_S$ carries the filtration induced by
$\deg T_1=1$ and $\deg F_h=h$. The filtered
polynomial extension gives
\begin{equation}\label{eq:filtered-sum}
(\ASp)_{\le D}=
\bigoplus_{0\le j\le\lfloor D/2\rfloor}T_2^j(B_S)_{\le D-2j},
\qquad \rhoSp(D)=\sum_{j\ge0}b_S(D-2j).
\tag{4.3}\end{equation}
Consequently $b_S(D)=\rhoSp(D)-\rhoSp(D-2)$.

\section{Intermediate growth and nonrationality}
We turn to the growth of the primitive filtration. The upper bound follows from the dimension of the symmetric-function ring. For the lower bound we need to construct many independent products. The construction also shows why the primitive spectrum detects divisibility by an increasing number of primes.

\begin{lemma}[A prime frequency close to the boundary]\label{lem:prime-frequency}
Write $g=\gcd(S)$, $S'=S/g$, and $M'=M/g$. There is a constant $C>0$ such that for every sufficiently large odd prime $\ell$ one can find $h=h(\ell)\ge3$ satisfying
\[
\ell\le hM'<\ell+C,\qquad h=\frac{\ell}{M'}+O_S(1),\qquad b_{h,g\ell}>0.
\]
\end{lemma}
\begin{proof}
Begin with $h$ copies of the largest normalized step $M'$. If we replace one by $\pm s'$, the frequency decreases by an element of
\[
E=\{M'-s',M'+s':s'\in S'\}.
\]
Let $\delta$ be the greatest common divisor of the positive
elements of $E$. It divides all $2s'$ and hence divides
$2\gcd(S')=2$. Thus $\delta=1$ or $2$. If $\delta=2$,
then every $M'-s'$ is even, so all $s'\in S'$ have the
same parity as $M'$. Since $\gcd(S')=1$, they must all
be odd; in particular $M'$ is odd. Consequently
$(M',\delta)=1$ in either case.

The additive semigroup generated by the positive defects contains every sufficiently large multiple of $\delta$. We choose $h$ so that $hM'-\ell$ belongs to a fixed bounded interval beyond the conductor and is divisible by $\delta$. Coprimality makes this possible with $hM'-\ell=O_S(1)$. This defect can be expressed using a bounded number of elements of $E$. Replacing that number of maximal steps gives a signed walk of length $h$ and displacement $\ell$. The Laurent coefficients are nonnegative, so $b_{h,g\ell}>0$.
\end{proof}

\begin{lemma}[Triangular separation of prime bits]\label{lem:prime-bits}
For every $k\ge1$ there exist $W_1,\ldots,W_k\in\ASp$, each of weight $O_S(i\log(i+1))$, such that the $2^k$ squarefree products
\[
\prod_{i\in I}W_i,\qquad I\subseteq\{1,\ldots,k\},
\]
are linearly independent as arithmetic functions.
\end{lemma}
\begin{proof}
Choose increasing odd primes $\ell_i>\max\{M,2M',g\}$
with gaps larger than the constant $C$ in
Lemma~\ref{lem:prime-frequency}. We may take every
$K$th sufficiently large prime for a fixed $K$.
Set $h_i=h(\ell_i)$ and arrange that
\[
\ell_i\le h_iM'<\ell_i+C<\ell_{i+1},\qquad h_iM'<2\ell_i.
\]
Take auxiliary primes $q,q'$ larger than $g$, $2M$, $h_kM'$, and every $\ell_i$. For each Boolean vector $x=(x_1,\ldots,x_k)\in\{0,1\}^k$, consider
\[
n_x=qq'\prod_{j=1}^k\ell_j^{x_j}.
\]
These numbers are admissible and squarefree, and every prime
divisor exceeds every element of $S$. Consequently
$c_{n_x}(s)=\mu(n_x)$ for each $s\in S$, whence
$T_1(n_x)=r\mu(n_x)$. Define
\[
W_i=\frac{T_1}{r}F_{h_i}.
\]
Its weight is at most $h_i+1$. For squarefree $n$ we have the elementary identity
\[
\mu(n)c_n(m)=\prod_{\substack{p\mid n\\p\mid m}}(1-p).
\]
Apply this to the Laurent expansions of the $F_{h_i}$.
Every nonzero frequency is a multiple of $g$, and its
normalized absolute value is at most $h_iM'$.
A later prime $\ell_j$, $j>i$, cannot divide it,
since $\ell_j>g$ and its normalized absolute value
is smaller than $\ell_{i+1}$.
Moreover, if $\ell_i$ divides one of these frequencies,
the normalized frequency must be $\pm\ell_i$, since
$h_iM'<2\ell_i$. The $R_2$ term in $F_{h_i}$
contributes no such frequency because $\ell_i>2M'$.
We therefore obtain, on the Boolean cube,
\[
W_i(n_x)=U_i(x_1,\ldots,x_{i-1})+\beta_i x_i,
\qquad \beta_i=-\ell_i(b_{h_i,g\ell_i}+b_{h_i,-g\ell_i})\ne0.
\]
Thus each $W_i$ introduces one new Boolean coordinate with a nonzero coefficient. By induction, the squarefree products form a basis of the functions on $\{0,1\}^k$; equivalently, their evaluation matrix is invertible. Finally,
\[
\sum_{i=1}^k(h_i+1)=O_S(k^2\log(k+1))
\]
by the prime number theorem. This gives the claimed weight bound and independence.
\end{proof}

\begin{theorem}[Universal primitive intermediate growth]\label{thm:prim-growth}
For every nonempty finite $S$, the estimates~\eqref{eq:prim-bounds} hold. The generating function $\sum_{D\ge0}\rhoSp(D)t^D$ and the Hilbert series of $\gr\ASp$ are both nonrational over $\Q(t)$.
\end{theorem}
\begin{proof}
Lemma~\ref{lem:prime-bits} gives $\rhoSp(D_k)\ge2^k$ for some $D_k=O_S(k^2\log(k+1))$. For large $D$, choose $k$ of order $\sqrt{D/\log D}$, with a sufficiently small constant so that $D_k\le D$. Since the filtered dimensions are nondecreasing, this proves the lower bound.

For the upper bound, $A_{S,\le D}^{\prim}$ is a quotient of $\Lambda_{\Q,\le D}$. The dimension of the latter is $\sum_{j\le D}p(j)$, where $p(j)$ is the partition function. The Hardy--Ramanujan estimate~\cite{HardyRamanujan} gives $\sum_{j\le D}p(j)=\exp(O(\sqrt D))$, which proves the upper bound.

The lower estimate grows faster than every polynomial, whereas the upper one is subexponential. A rational power series with nonnegative coefficients and subexponential coefficient growth can have only polynomially growing coefficients: this follows from its partial fractions, since no pole can lie strictly inside the unit disk. Therefore $\sum_D\rhoSp(D)t^D$ cannot be rational. The Hilbert series is obtained by multiplying it by $1-t$, so it is nonrational as well.
\end{proof}

\section{Quasirational rigidity and algebraic separation}\label{sec:rigidity}
We now prove Theorem~\ref{thm:rigidity-intro}. The difference between the two spectral algebras turns out to be stronger than the different growth estimates. We first show that a primitive spectral function is determined by any sufficiently long tail of its values.

\begin{lemma}[No finitely supported element]\label{lem:no-finite}
If $f\in\ASp$ vanishes at all sufficiently large integers, then $f$ vanishes at every admissible level $n>2M$.
\end{lemma}
\begin{proof}
Write $f=F(T_1,\ldots,T_H)$ for a polynomial $F\in\Q[X_1,\ldots,X_H]$. For every fixed admissible $d>2M$, define
\[
P_d(X)=F(a_1X+R_1(d),\ldots,a_HX+R_H(d))\in\Q[X].
\]
Choose two distinct primes $p,q>HM$, both coprime to $d$. Lemma~\ref{lem:insertion} gives
\[
f(dpq)=P_d\bigl(\varphi(d)(p-1)(q-1)\bigr).
\]
Keep $q$ fixed and vary $p$ over infinitely many large primes. By assumption, all the left-hand values are zero. The polynomial $P_d$ has infinitely many zeros and is therefore identically zero. Evaluating at $X=\varphi(d)$ gives $f(d)=0$. Since $d$ was arbitrary, the lemma follows.
\end{proof}

\begin{lemma}[Two-reflection rigidity]\label{lem:two-reflections}
If $H(X)\in\Q(X)$ is invariant under both $X\mapsto-X$ and $X\mapsto1-X$, then $H$ is constant.
\end{lemma}
\begin{proof}
Composing the two reflections gives translation by one, so $H(X+1)=H(X)$. A rational function in characteristic zero cannot have a nonzero period unless it is constant. Indeed, a finite pole would generate infinitely many poles by translation. With no finite poles the function is a polynomial, and a nonconstant polynomial cannot be periodic because of its leading term.
\end{proof}

\begin{theorem}[No nonconstant eventually quasirational primitive sequence]\label{thm:QR}
Let $f\in\ASp$. Suppose there exist $m\ge1$ and rational functions $Q_0,\ldots,Q_{m-1}\in\Q(X)$ such that for all sufficiently large admissible integers $n$,
\begin{equation}\label{eq:quasirational}
f(n)=Q_r(n)\qquad(n\equiv r\pmod m).
\tag{6.1}\end{equation}
Then $f$ is constant at every admissible level. Consequently $\ASp\cap\mathrm{QR}=\Q$.
\end{theorem}
\begin{proof}
Express $f$ as $F(T_1,\ldots,T_H)$ and consider the one-variable polynomial
\[
P(Y)=F(a_1Y,\ldots,a_HY)\in\Q[Y].
\]
The prime-square formula~\eqref{eq:insert-square} implies that for every sufficiently large prime $p$ coprime to $d$,
\begin{equation}\label{eq:prime-square-P}
f(dp^2)=P\bigl(\varphi(d)p(p-1)\bigr).
\tag{6.2}\end{equation}

For $m=1$, take any admissible $d_0$ and let $p$ run through all sufficiently large primes. Now suppose $m>1$. Infinitely many primes belong to one fixed invertible residue class $a$ modulo $m$: all but finitely many primes are coprime to $m$, and there are only finitely many invertible classes. Restrict $p$ to that infinite set. For each $r\pmod m$, choose $d_r>2M$ with
\[
d_r\equiv ra^{-2}\pmod m.
\]
In either case we have chosen $d_r$ for every residue class and an infinite set of primes $p$ for which $d_rp^2\equiv r\pmod m$.
Taking $p$ large enough also ensures $p\nmid d_r$. Equations~\eqref{eq:quasirational} and~\eqref{eq:prime-square-P} therefore give
\[
Q_r(d_rp^2)=P\bigl(\varphi(d_r)p(p-1)\bigr)
\]
for infinitely many primes $p$. Both sides are rational functions of $p$. Hence they agree identically:
\begin{equation}\label{eq:reflections-identity}
Q_r(d_rX^2)=P\bigl(\varphi(d_r)X(X-1)\bigr).
\tag{6.3}\end{equation}
The left side is unchanged by $X\mapsto-X$, while the right side is unchanged by $X\mapsto1-X$. Lemma~\ref{lem:two-reflections} shows that their common value must be constant. Since the substitution $X\mapsto\varphi(d_r)X(X-1)$ is nonconstant, $P$ itself is constant, say $P=c\in\Q$. Equation~\eqref{eq:reflections-identity} now gives $Q_r=c$ for each residue class $r$. Thus $f=c$ on all sufficiently large admissible integers. The function $f-c$ belongs to $\ASp$ and has finite support. Lemma~\ref{lem:no-finite} proves that it vanishes identically.
\end{proof}

\begin{corollary}[Trivial intersection in spite of prime-level agreement]\label{cor:intersection}
For every nonempty finite $S$,
\[
\ASg\cap\ASp=\Q.
\]
On the other hand, $G_h(p)=T_h(p)$ for every prime $p>2M$.
\end{corollary}
\begin{proof}
Every element of $\ASg$ is eventually polynomial by Proposition~\ref{prop:glob-image}, and so belongs to $\mathrm{QR}$. Theorem~\ref{thm:QR} gives the intersection. The equality on primes is immediate: every nonzero residue class modulo a prime is a unit.
\end{proof}
\begin{remark}
The global algebra contains every finitely supported sequence, while the primitive algebra contains none except zero. Thus the equality of all spectral moments at prime levels is compatible with the smallest possible intersection of the two algebras inside the common ring of rational arithmetic functions. This statement is algebraic; it does not assert any statistical independence.
\end{remark}

\section{Prime and rough-level quotients}
The additional complexity of the primitive algebra is found at composite moduli that have prime factors small in comparison with the filtration degree. When we restrict to large primes, or even to integers with only large prime factors, the picture becomes much simpler.

\begin{proposition}[Prime-tail quotient]\label{prop:prime-tail}
Modulo eventual equality on prime levels, the primitive evaluation image in degree at most $D$ has dimension
\[
\rho_{S,\mathrm{prime}}^{\prim}(D)=\left\lfloor\frac D2\right\rfloor+1.
\]
\end{proposition}
\begin{proof}
If $p>DM$ is prime, then $c_p(m)=-1$ for all nonzero $|m|\le DM$. Consequently $T_h(p)=a_hp-r^h=G_h(p)$ for every $h\le D$. Every evaluation of weight at most $D$ is therefore a polynomial in $p$ of degree at most $\lfloor D/2\rfloor$. These degrees all occur, because the powers of $T_2(p)=a_2p-r^2$ are independent polynomials.
\end{proof}

Call two arithmetic functions \emph{rough-tail equivalent} if, for some fixed bound, they agree on all sufficiently large integers whose least prime factor is above that bound. Let $A_{S,\mathrm{rough}}$ denote the corresponding quotient of $\ASp$.

\begin{theorem}[Exact rough-level algebra]\label{thm:rough}
There is a filtered algebra isomorphism
\begin{equation}\label{eq:rough-algebra}
A_{S,\mathrm{rough}}\simeq\Q[X,Y]/(Y^3-Y),\qquad \deg X=2,\quad\deg Y=1,
\tag{7.1}\end{equation}
where $X$ is represented by $\varphi(n)$ and $Y$ by $\mu(n)$. Its associated-graded Hilbert series equals $(1+t+t^2)/(1-t^2)$, and its filtered dimension grows linearly.
\end{theorem}
\begin{proof}
Suppose every prime dividing $n$ is larger than $hM$. Then $c_n(m)=\mu(n)$ for each nonzero frequency $|m|\le hM$. Formula~\eqref{eq:T-bridge} reduces to
\begin{equation}\label{eq:rough-moments}
T_h(n)=a_h\varphi(n)+(r^h-a_h)\mu(n).
\tag{7.2}\end{equation}
In particular,
\[
T_1=r\mu(n),\qquad T_2=a_2\varphi(n)+(r^2-a_2)\mu(n).
\]
Thus $Y=T_1/r$ has weight at most one, and
$X=(T_2-(r^2-a_2)Y)/a_2$ has weight at most two.
Conversely,~\eqref{eq:rough-moments} expresses every $T_h$
as a polynomial of weighted degree at most $h$ in
$X$ and $Y$ (with $a_1=0$ when $h=1$).
These maps therefore preserve the respective filtrations
in both directions. The relation $Y^3-Y=0$ follows from
$\mu(n)\in\{-1,0,1\}$.

There are no further polynomial relations. Indeed, for each $\varepsilon\in\{-1,0,1\}$, we can find rough integers with $\mu(n)=\varepsilon$ and infinitely many different values of $\varphi(n)$. Take large primes for $\varepsilon=-1$, squares of large primes for $\varepsilon=0$, and products of two distinct large primes for $\varepsilon=1$. Therefore, if $P(X,Y)$ vanishes on the rough tail, the three polynomial identities $P(X,\varepsilon)=0$ hold. This forces $P$ to be divisible by $Y(Y-1)(Y+1)$. The normal forms
\[
A(X)+YB(X)+Y^2C(X)
\]
give the vector-space basis
$\{X^j,YX^j,Y^2X^j:j\ge0\}$.
Their weights are $2j,2j+1,2j+2$,
which proves the associated-graded Hilbert
series $(1+t+t^2)/(1-t^2)$.
\end{proof}

\section{An exact finite-matrix formula}
\begin{theorem}[Finite evaluation-matrix rank]\label{thm:finite-matrix}
Let $S$ be any nonempty finite jump set and fix $D\ge1$. Put
\[
B=DM,\qquad L_B=\lcm(1,\ldots,B),\qquad
K_B=L_B\prod_{p\le B}p.
\]
Index the columns by all partitions
$\lambda=(\lambda_1,\ldots,\lambda_t)$ with
$0\le|\lambda|\le D$, and the rows by triples
$(d,\varepsilon,j)$ with $d\mid K_B$,
$\varepsilon\in\{-1,0,1\}$, and
$0\le j\le\lfloor D/2\rfloor$. Define the finite
matrix over $\Q$ by
\[
M_{S,D}[(d,\varepsilon,j),\lambda]
=[X^j]\prod_{i=1}^t\bigl(a_{\lambda_i}X+\varepsilon R_{\lambda_i}(d)\bigr),
\]
where the empty partition gives the constant polynomial $1$.
Then
\begin{equation}\label{eq:finite-matrix-rank}
\rhoSp(D)=\rank_{\Q}M_{S,D}.
\tag{8.3}\end{equation}
\end{theorem}
\begin{proof}
For each column partition $\lambda$, write
$T_\lambda=\prod_iT_{\lambda_i}$, with
$T_\varnothing=1$. These monomials span
$(\ASp)_{\le D}$. We show that a linear combination
\[
f=\sum_{|\lambda|\le D}c_\lambda T_\lambda
\]
vanishes on every admissible level if and only if
the column vector $(c_\lambda)_\lambda$ belongs to
$\ker M_{S,D}$.

For $d\mid K_B$ and $\varepsilon\in\{-1,0,1\}$, set
\[
P_{d,\varepsilon}(X)=
\sum_{|\lambda|\le D}c_\lambda
\prod_i\bigl(a_{\lambda_i}X+
\varepsilon R_{\lambda_i}(d)\bigr).
\]
This polynomial has degree at most $\lfloor D/2\rfloor$,
because $a_1=0$ and every factor with a nonzero
coefficient of $X$ has weight at least two.
Consequently $M_{S,D}(c_\lambda)_\lambda=0$
is equivalent to
$P_{d,\varepsilon}=0$ for every pair $(d,\varepsilon)$.

Suppose first that $f$ is the zero arithmetic function.
Fix $d\mid K_B$. There are infinitely many integers
$q$ coprime to $d$, with all prime factors greater than
$B$, such that $\mu(q)=\varepsilon$ and
$\varphi(q)$ takes infinitely many distinct values.
Indeed, take a varying large prime if $\varepsilon=-1$,
the square of a varying large prime if $\varepsilon=0$,
or the product of a fixed large prime and a second
varying large prime if $\varepsilon=1$.
For every $h\le D$, the multiplicativity of Ramanujan
sums in coprime moduli gives
\[
R_h(dq)=\mu(q)R_h(d)=\varepsilon R_h(d),
\qquad
\varphi(dq)=\varphi(d)\varphi(q),
\]
because no prime factor of $q$ divides a nonzero
frequency of absolute value at most $hM\le B$.
For all sufficiently large such $q$, we have
\[
0=f(dq)=
P_{d,\varepsilon}\bigl(\varphi(d)\varphi(q)\bigr).
\]
Thus $P_{d,\varepsilon}$ has infinitely many
distinct roots and is identically zero.
This proves $(c_\lambda)_\lambda\in\ker M_{S,D}$.

Conversely, assume every $P_{d,\varepsilon}$ vanishes.
Fix any admissible $n>2M$ and factor it as
$n=n_{\le B}n_{>B}$, where $n_{\le B}$ is supported
on primes at most $B$ and $n_{>B}$ has all prime
factors greater than $B$.
If $n_{\le B}\mid K_B$, set
$d=n_{\le B}$ and $\varepsilon=\mu(n_{>B})$.
The same multiplicativity argument gives
\[
f(n)=P_{d,\varepsilon}\bigl(\varphi(n)\bigr)=0.
\]
Otherwise some prime $p\le B$ occurs in $n$ to an
exponent greater than
\[
v_p(K_B)=\lfloor\log_p B\rfloor+1.
\]
Every nonzero frequency $m$ appearing in any $R_h$
with $h\le D$ satisfies $|m|\le B$, and hence
$v_p(m)\le\lfloor\log_p B\rfloor$.
It follows that
$v_p(n)\ge v_p(m)+2$, so the prime-power
Ramanujan sum $c_{p^{v_p(n)}}(m)$ is zero.
By multiplicativity, $c_n(m)=0$, whence
$R_h(n)=0$ for every $h\le D$.
Therefore
\[
f(n)=P_{1,0}\bigl(\varphi(n)\bigr)=0.
\]
We have proved equality of the kernels of the
arithmetic evaluation map and $M_{S,D}$.
Equality of their ranks proves~\eqref{eq:finite-matrix-rank}.
\end{proof}

\begin{remark}[The zero-correction row]
The second case in the proof is necessary.
For example, if $S=\{1\}$ and $D=2$, then
$B=2$, $K_B=4$, and $R_2(n)=c_n(2)/2$.
Thus $R_2(8)=0$ but $R_2(4)=-1$.
One cannot obtain the correction at $n=8$
simply by truncating the power of $2$ in $n$
to the corresponding power in $K_B$.
Instead, all relevant corrections vanish and
the row $(d,\varepsilon)=(1,0)$ applies.
\end{remark}

\section{Discussion and open problems}
The two natural evaluations of the same family of circulant spectra
are algebraically very different. The full-spectrum functions are
eventually polynomial; the primitive moments are sensitive to
increasingly many prime divisibility conditions. This sensitivity
produces intermediate filtered growth and a nonrational Hilbert
series. Even though the two evaluations agree at every prime order,
their algebras have only the constants in common.

The bounds for $\log\rho_S^{\prim}(D)$ leave a substantial gap
between $\sqrt{D/\log D}$ and $\sqrt D$. Determining the correct
growth scale, and describing the locally finite algebra $B_S$
intrinsically rather than through finite evaluation matrices, remain
open problems. The intermediate growth refers specifically to
the power-sum filtration: it is distinct from the linear word growth
of every finitely generated subalgebra.

\section*{AI-assisted preparation}

OpenAI ChatGPT (GPT-6) was used as a supporting tool
in mathematical discussions, critical examination
of proof details, preparation and revision of draft
passages, and English-language and LaTeX editing.

In particular, AI-assisted discussion helped
identify a gap in an earlier argument concerning
the finite-matrix rank formula.

The AI system is not an author. The named authors
retain full responsibility for the mathematical
content, originality, accuracy of the references,
and final submitted version of the manuscript.

\end{document}